\documentclass[11pt]{amsart}

\usepackage{amscd}
\usepackage{amsmath}
\usepackage{graphicx}
\usepackage{amsfonts}
\usepackage{amssymb}
\newcommand{\tr}{{\rm tr}\,}

\begin{document}
\title[Maps preserving peripheral spectrum]
{Maps preserving peripheral spectrum of \\generalized Jordan products
of operators}

\author{Wen Zhang}
\address[Wen Zhang ]{Department of
 Mathematics, Shanxi University, Taiyuan,
Shanxi,  030006, P. R. China}\email{wenzhang1314gw@163.com}

\author{Jinchuan Hou}
\address[Jinchuan Hou]{Department of Mathematics, Taiyuan
University of Technology, Taiyuan 030024, P. R. of China;
Department of
 Mathematics, Shanxi University, Taiyuan,
Shanxi,  030006, P. R. China} \email{houjinchuan@tyut.edu.cn}

\author{Xiaofei Qi}
\address[Xiaofei Qi]{Department of
 Mathematics, Shanxi University, Taiyuan,
Shanxi,  030006, P. R. China}\email{qixf1980@126.com}

\thanks{This work is partially  supported by National Natural Science Foundation
of China (No.11171249, 11101250, 11271217).}

\begin{abstract}
Let $X_1$ and $X_2$ be complex Banach spaces with dimension at least
three, $\mathcal{A}_1$ and $\mathcal{A}_2$ be standard operator
algebras on $X_1$ and $X_2$, respectively. For $k\geq2$, let
$(i_1,...,i_m)$ be a sequence with terms chosen from
$\{1,\ldots,k\}$ and assume that at least one of the terms in
$(i_1,\ldots,i_m)$ appears exactly once. Define the generalized
Jordan product $T_1\circ T_2\circ\cdots\circ T_k=T_{i_1}
T_{i_2}\cdots T_{i_m}+T_{i_m}\cdots T_{i_2} T_{i_1}$ on elements in
$\mathcal{A}_i$. This includes the usual Jordan product
$A_1A_2+A_2A_1$, and the Jordan triple $A_1A_2A_3+A_3A_2A_1$. Let
$\Phi:\mathcal{A}\rightarrow\mathcal{B}$ be a map  with range
containing all operators of rank at most three. It is shown that
$\Phi$ satisfies that
$\sigma_\pi(\Phi(A_1)\circ\cdots\circ\Phi(A_k))=\sigma_\pi(A_1\circ\cdots\circ
A_k)$ for all $A_1, \ldots, A_k$, where $\sigma_\pi(A)$ stands for
the peripheral spectrum of $A$, if and only if   $\Phi$ is a Jordan
isomorphism multiplied by an $m$th root of unity.

\end{abstract}
\maketitle

\section{Introduction}

There has been considerable interest in studying spectrum preserving
maps on operator algebras in connection to the Kaplansky¡¯s problem
on characterization of linear maps between Banach algebras
preserving invertibility; see \cite{A1,AM,JS,K1,S1}. Early study
focus on linear maps, additive maps, or multiplicative maps; see,
e.g., \cite{LT2}. Moreover, spectrum preserving maps on Banach
algebras which are not assumed to be linear are studied by several
authors (see \cite{LLT,LT1,M1,PS,R1,RR1,RR2}). In \cite{M1},
Moln$\acute{a}$r characterized surjective maps $\Phi$ on bounded
linear operators acting on a Hilbert space preserving the spectrum
of the product of operators, i.e., $AB$ and $\Phi(A)\Phi(B)$ always
have the same spectrum. Hou, Li and Wong \cite{HLW1,HLW2} studied
respectively further the maps $\Phi$ between certain operator
algebras preserving the spectrum of a generalized product $T_1*
T_2*\cdots* T_k$ and a generalized Jordan product $T_1\circ
T_2\circ\cdots\circ T_k$ of low rank operators.

{\bf Definition 1.1.} {\it Fix a positive integer $k$ and a finite
sequence $(i_1,i_2,\ldots,i_m)$ such that
$\{i_1,i_2,\ldots,i_m\}=\{1,2,\ldots,k\}$ and there is an $i_p$ not
equal to $i_q$ for all other $q$. For operators $T_1,\ldots,T_k$,
the operators
$$T_1* T_2*\cdots* T_k= T_{i_1}T_{i_2}\cdots
T_{i_m} \eqno(1.1)$$ and
  $$T_1\circ T_2\circ\cdots\circ T_k= T_{i_1}T_{i_2}\cdots
T_{i_m}+T_{i_m}\cdots T_{i_2} T_{i_1} \eqno(1.2)$$ are respectively
called generalized product and generalized Jordan product of
$T_1,\ldots,T_k$.}

Evidently, the generalized product $T_1*\cdots*T_k$ (the generalized
Jordan product $T_1\circ\cdots\circ T_k$) covers the usual product
$T_1T_2$ and the Jordan semi-triple product $T_1T_2T_1$ (the Jordan
product $T_1T_2 + T_2T_1$ and the Jordan triple product $T_1T_2T_3 +
T_3T_2T_1$), etc.. In \cite{HLW1} (resp. \cite{HLW2}) it was shown
that, if $\Phi$ is a  map between standard operator algebras with
the range  containing all operators of rank at most three, then, for
all operators $T_1, T_2,\cdots,T_k$ of low rank, the spectra of
$T_1* T_2*\cdots* T_k$ (resp., of $T_1\circ T_2\circ\cdots\circ
T_k$)  and of $\Phi(T_1)*\Phi(T_2)*\cdots* \Phi(T_k)$ (resp., of
$\Phi(T_1)\circ\Phi(T_2)\circ\cdots\circ \Phi(T_k)$) are equal if
and only $\Phi$ is a Jordan isomorphism multiplied by a $m$th root
of the unit.

Let $\mathcal{B}(X)$ be the Banach algebra of all bounded linear
operators on a complex Banach space $X$. Denote by $\sigma(T)$ and
$r(T)$ the spectrum and the spectral radius of $T\in \mathcal{B}(X)$,
respectively. The peripheral spectrum of $T$ is defined by
$$\sigma_\pi(T)=\{z\in\sigma (T):|z|=r(T)\}.$$ Since $\sigma(T)$ is
compact, $\sigma_\pi(T)$ is a well-defined non-empty set and is an
important spectral function. In \cite{TL}, Tonev and Luttman studied
maps preserving peripheral spectrum of the usual operator products
on standard operator algebras. They studied also the corresponding
problems in uniform algebras  (see \cite{LLT,LT1}).  Later Takeshi
and Dai \cite{TD} generalized the result in \cite{LT1} and
characterized surjective maps $\phi$ and $\psi$ satisfying
$\sigma_\pi(\phi(T)\psi(S))=\sigma_\pi(TS)$ on standard operator
algebras. In \cite{ZH1} the  maps preserving peripheral spectrum of
Jordan semi-triple products of operators is characterized,  and
then, the maps preserving peripheral spectrum of generalized
products of operators is also characterized in \cite{ZH2}. The maps
preserving peripheral spectrum of Jordan products  $AB+BA$ of
operators on standard operator algebras are  characterized by Cui
and Li
 in \cite{CL}.

In this paper, we continue the study and characterize the maps
 preserving the peripheral spectrum of generalized Jordan products of
 operators between the standard operator algebras on complex Banach spaces. As expected, such maps
 are
 Jordan isomorphisms multiplied by a suitable root of the unit 1.
 However, the situation of generalized Jordan product is much more
 complicated than the Jordan product case.  We
can not use the similar  technique for the case of Jordan product
$AB+BA$ discussed in \cite{CL} to solve our problem. One of the
reasons is that the   root of $\Phi(I)$ may not be obtained.

In the following, let $X_i$ be a complex Banach space with dimension
at least three  and $\mathcal{A}_i$ the standard operator algebra on
$X_i$, i.e., $\mathcal{A}_i$ contains all continuous finite rank
operators on $X_i$,  $i = 1, 2$. Note that a Jordan isomorphism
$\Phi: \mathcal{A}_1\rightarrow\mathcal{A}_2$ is either an inner
automorphism or anti-automorphism. In this case, it is obvious that
$\sigma_\pi(\Phi(A_1)\circ
\cdots\circ\Phi(A_k))=\sigma_\pi(A_1\circ\cdots\circ A_k)$ holds for
all $A_1, \ldots,A_k$. The main result of this paper is to show that
the converse is also true.

{\bf Theorem 1.2.} {\it Let ${\mathcal A}_i$ be a standard operator
algebra on a complex Banach space $X_i$ with $\dim X_i\geq 3$,
$i=1,2$. Consider the product $T_1\circ\cdots\circ T_k$ defined in
Definition 1.1. Let $\Phi:\mathcal{A}_1\rightarrow\mathcal{A}_2$ be
a map with the range   containing all operators in $\mathcal{A}_2$
of rank at most three. Then $\Phi$ satisfies
$$\sigma_\pi(\Phi(A_1)\circ\cdots\circ\Phi(A_k))=\sigma_\pi
(A_1\circ\cdots\circ A_k) \eqno(1.3)$$ for any $A_1,\ldots,
A_k\in{\mathcal A}_1$  if and only if  one of the following
conditions holds.}

(1) {\it There exist a scalar $\lambda\in\mathbb{C}$ with $\lambda^m=1$
and an invertible operator $T\in\mathcal{B}(X_1,X_2)$
such that $\Phi(A)=\lambda TAT^{-1}$ for all $A\in\mathcal{A}_1$.}

(2) {\it The spaces $X_1$ and $X_2$ are reflexive, and there exist a
scalar $\lambda\in\mathbb{C}$ with $\lambda^m=1$ and an invertible
operator $T\in\mathcal{B}(X^*_1,X_2)$ such that $\Phi(A)=\lambda
TA^*T^{-1}$ for all $A\in\mathcal{A}_1$.}

The ``if" part is clear. So we need only to show the ``only if" part
of Theorem 1.2.  This will be done  in Section 3.

\section{Characterizations of rank one operators}

To prove  Theorem 1.2, it is important to characterize   rank one
operators in terms of the peripheral spectrum and generalized Jordan
products.

{\bf Lemma 2.1.} {\it  Let $X$ be a complex Banach space  with $\dim
X\geq 3$. Assume that  $A\in{\mathcal B}(X)$ is a nonzero operator,
$r$ and $s$ are integers with $s>r>0$.   Then the following
conditions are equivalent.}

(1){ \it $A$ has rank one.}

(2) {\it $\sigma_\pi(B^rAB^s+B^sAB^r)$ has at most two elements for
any $B\in{\mathcal B}(X)$.}

(3){ \it There does not exist an operator $B\in{\mathcal B}(X)$ with
rank at most three such that $B^rAB^s+B^sAB^r$ has rank three and
$\sigma_\pi(B^rAB^s+B^sAB^r)$ has at most two elements.}

{\bf Proof}. The approach is similar to \cite{HLW2} but more
complicated. We give details for reader's convenience.

The implications (1)$\Rightarrow$(2)$\Rightarrow$(3) are clear.

To prove $(3)\Rightarrow(1)$, suppose (3) holds but (1) is not
true, i.e., $A$ has rank at least two.

If $A$ has rank at least 3, then there are $x_1, x_2, x_3\in X$ such
that $\{Ax_1,Ax_2,Ax_3\}$ is linearly independent. Consider the
operator matrix of $A$ on ${\rm span}\{x_1, x_2,
x_3,Ax_1,Ax_2,Ax_3\}$ and its complement:
$$\left(\begin{array}{cc} A_{11}&A_{12}\\A_{21}&A_{22}\end{array}\right).$$
Then $A_{11}\in M_n$ with $3\leq n\leq6$. By \cite{HLW1}, there is a
nonsingular $U$ on  ${\rm span}\{x_1, x_2, x_3,Ax_1,$ $Ax_2,Ax_3\}$
such that $U^{-1}A_{11}U$ has an invertible 3-by-3 leading
submatrix. We may further assume that the 3-by-3 matrix is in
triangular form with nonzero diagonal entries $a_1, a_2, a_3$. Now
let $B\in\mathcal{A}$ has operator matrix
$$\left(\begin{array}{cc} B_{11}& 0\\0&0\end{array}\right),$$
where $UB_{11}U^{-1}$ = diag $(1, b_2, b_3)\oplus 0_{n-3}$
 with $B_{11}$ using the same basis as that of $A_{11}$ and $b_2, b_3$ being chosen
such that $a_1, a_2b^{r+s}_2 , a_3b^{r+s}_3$ are three distinct nonzero numbers with
$|a_1|=| a_2b^{r+s}_2|=|a_3b^{r+s}_3|$. It follows
that $B^rAB^s+B^sAB^r$ has rank 3 and $\sigma_\pi(B^rAB^s + B^sAB^r)$ has three different
points.

Next, suppose $A$ has rank 2. Choosing a suitable space
decomposition of $X$, we may assume that $A$ has operator matrix
$A_1\oplus 0$, where $A_1$ has one of the following forms.
$$\rm{(i)}\left(\begin{array}{ccc} a&0&b\\0&0&0\\0&0&c\end{array}\right),
\rm{(ii)}\left(\begin{array}{ccc} a&0&0\\0&0&1\\0&0&0\end{array}\right),
\rm{(iii)}\left(\begin{array}{ccc} 0&1&0\\0&0&1\\0&0&0\end{array}\right),
\rm{(iv)}\left(\begin{array}{cc} 0_2&I_2\\0_2&0_2\end{array}\right).$$

If (i) holds, set $\theta=\frac{\pi}{2(r+s)}$. Then
$\cos(r+s)\theta=0$, $\sin(r+s)\theta=1$ and
$\sin(r-s)\theta\neq0$. Let $d$ be a number such that
$|a\sin(r-s)\theta|=|2cd^{r+s}|$ but $2cd^{r+s}\neq\pm
ia\sin(r-s)\theta$. Let $B\in \mathcal{A}$ be represented by the
operator matrix
$$\left(\begin{array}{ccc}
\cos\theta&-\sin\theta&0\\
\sin\theta&\cos\theta&0\\
0&0&d\end{array}\right)\oplus0.$$
Then $B^rAB^s+B^sAB^r$ has operator matrix
$$\left(\begin{array}{ccc}
a\cos(r-s)\theta&-a&*\\a&-\cos(r-s)\theta&*\\
0&0&2cd^{r+s}\end{array}\right)\oplus0,$$
which has rank 3 and $\sigma_\pi(B^rAB^s + B^sAB^r)=\{2cd^{r+s},ia\sin(r-s)\theta,
-ia\sin(r-s)\theta\}$.

Suppose (ii) holds, set $\theta=\frac{\pi}{r+s}$. Then
$\cos(r+s)\theta=-1$, $\sin(r+s)\theta=0$ and
$\sin(r-s)\theta\neq0$. Let $d$ be a number such that
$|\sin(r-s)\theta|=|2ad^{r+s}|$ but $2ad^{r+s}\neq\pm
i\sin(r-s)\theta$. Constructing $B$ by the operator matrix
$$\left(\begin{array}{ccc}
d&0&0\\
0&\cos\theta&-\sin\theta\\
0&\sin\theta&\cos\theta\end{array}\right)\oplus0,$$
$B^rAB^s+B^sAB^r$ has operator matrix
$$\left(\begin{array}{ccc}
2ad^{r+s}&0&0\\
0&0&\cos(r-s)\theta-1\\
0&\cos(r-s)\theta+1&0\end{array}\right)\oplus0,$$
which has rank 3 and $\sigma_\pi(B^rAB^s + B^sAB^r)=\{2ad^{r+s},i\sin(r-s)\theta,
-i\sin(r-s)\theta\}$.

Suppose (iii) holds. First, assume that $s = 2r$. Let $B$ be such that $B^r$
has operator matrix
$$\left(\begin{array}{ccc}
0&1&0\\
0&0&1\\
1&0&0\end{array}\right)\oplus0.$$
Then $B^rAB^s+B^sAB^r$ has operator matrix
$$\left(\begin{array}{ccc}
0&1&0\\
0&0&1\\
2&0&0\end{array}\right)\oplus0,$$
which has rank 3 and $\sigma_\pi(B^rAB^s + B^sAB^r)=\{2^{1/3},2^{1/3}e^{i2\pi/3},2^{1/3}e^{i4\pi/3}\}.$

Next, suppose $s/r \neq 2$. Then $s>2$ and $2r/s$ is not an integer. Let $\theta_1=2\pi/s,\theta_2=4\pi/s.$
Then $1, e^{ir\theta_1}, e^{ir\theta_2}$ are distinct because $e^{i4\pi r/s}=e^{i2\pi(2r/s)}\neq1$
and $e^{ir\theta_1}=e^{ir\theta_2}/e^{ir\theta_1}=e^{i2\pi r/s}\neq1$.
 Thus, there exists an invertible $S\in M_3$ such that
$$\left(\begin{array}{ccc}
1&0&0\\
0&e^{ir\theta_1}&0\\
0&0&e^{ir\theta_2}\end{array}\right)
=S^{-1}\left(\begin{array}{ccc}
1&\alpha&0\\
0&e^{ir\theta_1}&0\\
m&1&e^{ir\theta_2}\end{array}\right)S.$$
Let $B$ have operator matrix
$$S\left(\begin{array}{ccc}
1&0&0\\
0&e^{i\theta_1}&0\\
0&0&e^{i\theta_2}\end{array}\right)S^{-1}\oplus0.$$
The operator matrix $B^s = I_3\oplus0$ and the operator matrix of $B^r$ has the form
$$S\left(\begin{array}{ccc}
1&0&0\\
0&e^{ir\theta_1}&0\\
0&0&e^{ir\theta_2}\end{array}\right)S^{-1}\oplus0
=\left(\begin{array}{ccc}
1&\alpha&0\\
0&e^{ir\theta_1}&0\\
m&1&e^{ir\theta_2}\end{array}\right)\oplus0.$$
Then $B^rAB^s + B^sAB^r = AB^r + B^rA$ has operator matrix
$$\left(\begin{array}{ccc}
0&1+e^{ir\theta_1}&\alpha\\
m&1&e^{ir\theta_1}+e^{ir\theta_2}\\
0&m&1\end{array}\right)\oplus0,$$
which has rank 3. It follows from $\cos\frac{2r\pi}{s}+1\neq0$ that $2e^{ir\theta_1}+e^{ir\theta_2}+1=2
(\cos\frac{2r\pi}{s}+i\sin\frac{2r\pi}{s})
+(\cos\frac{4r\pi}{s}+i\sin\frac{4r\pi}{s})+1=2\cos\frac{2r\pi}{s}(\cos\frac{2r\pi}{s}+1)
+i2\sin\frac{2r\pi}{s}(\cos\frac{2r\pi}{s}+1)\neq0$. So
let $m=\frac{-3}{2e^{ir\theta_1}+e^{ir\theta_2}+1}$ and $\alpha=\frac{(2e^{ir\theta_1}+e^{ir\theta_2}+1)
(5+13e^{ir\theta_1}+8e^{ir\theta_2})}{9}$, then $\sigma_\pi(B^rAB^s + B^sAB^r) =\sigma_\pi( AB^r + B^rA)
=\{2i,-2i,2\}$.

If (iv) holds, then $X$ has dimension at least 4. We may use a different decomposition of $X$
and assume that $A$ has operator matrix
$$\left(\begin{array}{cc} 0&1\\0&0\end{array}\right)\oplus
\left(\begin{array}{cc} 1&1\\-1&-1\end{array}\right)\oplus0.$$ Let
$\theta=\frac{\pi}{r+s}$. Then $\cos(r+s)\theta=-1$,
$\sin(r+s)\theta=0$ and $\sin(r-s)\theta\neq0$. Let $d$ be such
that $|\sin(r-s)\theta|=|d^{r+s}|$ but $d^{r+s}\neq\pm
i\sin(r-s)\theta$. Let $B\in \mathcal{A}$ be represented by the
operator matrix
$$\left(\begin{array}{ccc}
\cos\theta&-\sin\theta&0\\
\sin\theta&\cos\theta&0\\
0&0&d\end{array}\right)\oplus0.$$
Then $B^rAB^s+B^sAB^r$ has operator matrix
$$\left(\begin{array}{ccc}
0&\cos(r-s)\theta+\cos(r+s)\theta&0\\
\cos(r-s)\theta-\cos(r+s)\theta&0&0\\
0&0&d^{r+s}\end{array}\right)\oplus0,$$
which has rank 3 and $\sigma_\pi(B^rAB^s + B^sAB^r)=\{d^{r+s},i\sin(r-s)\theta,
-i\sin(r-s)\theta\}$.\hfill$\Box$

{\bf Lemma 2.2.} {\it  Suppose $s$ is a positive integer. Let $X$ be
a complex Banach space with $\dim X\geq 3$. Let $A\in\mathcal{B}(X)$
be such that $A^2\neq0$. Then the following conditions are
equivalent.}

(1){ \it $A$ has rank one.}

(2) {\it $\sigma_\pi(AB^s+B^sA)$ has at most two elements for any
$B$ in ${\mathcal B}(X)$.}

(3){ \it $\sigma_\pi(AB^s+B^sA)$ has at most two elements for any
$B$ whenever ${\rm rank}B\leq3$ and ${\rm rank}(AB^s+B^sA)\leq3$.}

{\bf Proof}. The implications (1)$\Rightarrow$(2)$\Rightarrow$(3) are
clear.

Suppose (3) holds but (1) is not true, i.e., $A$ has rank at least
two such that $A^2\neq0$.

First suppose $A$ has rank 2. Since $A^2\neq0$, choosing a
suitable space decomposition of $X$, we may assume that $A$ has
operator matrix $A_1\oplus 0$, where $A_1$ has one of the
following forms.
$$\rm{(i)}\left(\begin{array}{ccc} a&0&b\\0&0&0\\0&0&c\end{array}\right),
\rm{(ii)}\left(\begin{array}{ccc}
a&0&0\\0&0&1\\0&0&0\end{array}\right),
\rm{(iii)}\left(\begin{array}{ccc}
0&1&0\\0&0&1\\0&0&0\end{array}\right).
$$

 If (i) holds, set $\theta=\frac{\pi}{2s}$.
Then $\cos s\theta=0$, $\sin s\theta=1$. Let $d$ be such that
$|a|=|2cd^s|$ but $2cd^s\neq\pm ia$. Let $B\in \mathcal{A}$ be
represented by the operator matrix
$$\left(\begin{array}{ccc}
\cos\theta&-\sin\theta&0\\
\sin\theta&\cos\theta&0\\
0&0&d\end{array}\right)\oplus0.$$
Then $AB^s+B^sA$ has operator matrix
$$\left(\begin{array}{ccc}
0&-a&*\\a&0&*\\
0&0&2cd^s\end{array}\right)\oplus0,$$ which implies that
$\sigma_\pi(AB^s + B^sA)=\{2cd^s,ia,-ia\}$.

Suppose (ii) holds. Since the matrix$$C=\left(\begin{array}{ccc}
1/2&a&0\\
0&0&0\\
-1/2&0&-2\end{array}\right)$$
is similar to a matrix with distinct eigenvalues 0, 1/2, -2,
there exists an operator $B$ of rank 2 such that the operator matrix
of $B^s$ equals $C\oplus0$. It follows that the operator matrix of
$AB^s + B^sA$ is
$$\left(\begin{array}{ccc}
a&a^2&a\\
-1/2&0&-2\\
-a/2&0&0\end{array}\right)\oplus0,$$
and $\sigma_\pi(AB^s + B^sA)=\{a,ia,-ia\}$.

Suppose (iii) holds. Since the matrix$$C=\left(\begin{array}{ccc}
0&1&0\\
0&0&1\\
1&0&0\end{array}\right)$$
has distinct eigenvalues $1,e^{i2\pi/3},e^{i4\pi/3}$, there exists an operator $B$
of rank 3 such that the operator matrix of $B^s$ equals $C\oplus0$.
Then $AB^s + B^sA$ has operator matrix $$\left(\begin{array}{ccc}
0&0&2\\
1&0&0\\
0&1&0\end{array}\right)\oplus0,$$
and $\sigma_\pi(AB^s + B^sA)=\{2^{1/3},2^{1/3}e^{i2\pi/3},2^{1/3}e^{i4\pi/3}\}$.

Now, suppose $A$ has rank at least 3. Since $A^2\neq0$, there is $x\in X$ such that
$A^2x\neq0$. We consider the following two cases.

{\bf Case 1.} There is $x\in X$ such that $[x,Ax,A^2x]$ has dimension 3.

Decompose $X$ into $[x,Ax,A^2x]$ and its complement.
The operator matrix of $A$ has the form
$$\left(\begin{array}{cccc}
0&0&c_1&*\\
1&0&c_2&*\\
0&1&c_3&*\\
0&0&*&*\end{array}\right).$$

{\bf Subcase 1.} $c_1\neq0,c_2=0$.

Since the matrix$$C=\left(\begin{array}{ccc}
1&0&0\\
0&2&0\\
0&0&0\end{array}\right)$$
has distinct eigenvalues $1,2,0$, there exists an operator $B$
of rank 2 such that the operator matrix of $B^s$ equals $C\oplus0$.
Then $AB^s + B^sA$ has operator matrix $$\left(\begin{array}{cccc}
0&0&c_1&*\\
3&0&0&*\\
0&2&0&0\\
0&0&0&0\end{array}\right).$$ Writing $6c_1=re^{i\theta}$, we have
$\sigma_\pi(AB^s + B^sA)=\{r^{1/3}e^{i\theta/3},
r^{1/3}e^{i(\theta+2\pi)/3},r^{1/3}e^{i(\theta+4\pi)/3}\}$.

{\bf Subcase 2.} $c_1\neq0,c_2\neq0,c_3=0$.

Let $\alpha\neq0$ such that $\alpha^2c_2+1=0$,
Since the matrix$$C=\left(\begin{array}{ccc}
\alpha&1&0\\
0&-\alpha&-1\\
0&0&0\end{array}\right)$$
has distinct eigenvalues $\alpha,-\alpha,0$, there exists an operator $B$
of rank 2 such that the operator matrix of $B^s$ equals $C\oplus0$.
Then $AB^s + B^sA$ has operator matrix $$\left(\begin{array}{cccc}
1&0&\alpha c_1+c_2&*\\
0&0&-\alpha c_2&*\\
0&-\alpha&-1&0\\
0&0&0&0\end{array}\right),$$
and $\sigma_\pi(AB^s + B^sA)=\{1,e^{i2\pi/3},e^{i4\pi/3}\}$.

{\bf Subcase 3.} $c_1\neq0,c_2\neq0,c_3\neq0$.

Let $\alpha\neq0$ such that $3\alpha c_1+4c_2=0$. Since the
matrix$$C=\left(\begin{array}{ccc}
1&0&0\\
\alpha&2&0\\
0&0&0\end{array}\right)$$
has distinct eigenvalues $1,2,0$, there exists an operator $B$
of rank 2 such that the operator matrix of $B^s$ equals $C\oplus0$.
Then $AB^s + B^sA$ has operator matrix $$\left(\begin{array}{cccc}
0&0&c_1&*\\
3&0&2c_2+\alpha c_1&*\\
\alpha&2&0&0\\
0&0&0&0\end{array}\right).$$ With $6c_1=re^{i\theta}$, we see that
$\sigma_\pi(AB^s + B^sA)=\{r^{1/3}e^{i\theta/3},
r^{1/3}e^{i(\theta+2\pi)/3},r^{1/3}e^{i(\theta+4\pi)/3}\}$.

{\bf Subcase 4.} $c_1=0,c_2\neq0,c_3\neq0$.

Choose $\alpha\neq0$ such that $\alpha^2c_2+\alpha c_3+1=0$. Since
the matrix $$C=\left(\begin{array}{ccc}
\alpha&1&0\\
0&-\alpha&-1\\
0&0&0\end{array}\right)$$
has distinct eigenvalues $\alpha,-\alpha,0$, there exists an operator $B$
of rank 2 such that the operator matrix of $B^s$ equals $C\oplus0$.
Then $AB^s + B^sA$ has operator matrix $$\left(\begin{array}{cccc}
1&0&c_2&*\\
0&0&-\alpha c_2-c_3&*\\
0&-\alpha&-1&0\\
0&0&0&0\end{array}\right),$$ and thus $\sigma_\pi(AB^s +
B^sA)=\{1,e^{i2\pi/3},e^{i4\pi/3}\}$.

{\bf Subcase 5.} $c_1=0,c_2=0,c_3\neq0$.

Let $\alpha\neq0$ such that $\alpha c_3+1=0$,
Since the matrix$$C=\left(\begin{array}{ccc}
\alpha&1&0\\
0&-\alpha&-1\\
0&0&0\end{array}\right)$$
has distinct eigenvalues $\alpha,-\alpha,0$, there exists an operator $B$
of rank 2 such that the operator matrix of $B^s$ equals $C\oplus0$.
Then $AB^s + B^sA$ has operator matrix $$\left(\begin{array}{cccc}
1&0&0&*\\
0&0&c_3&*\\
0&-\alpha&-1&0\\
0&0&0&0\end{array}\right),$$ and thus $\sigma_\pi(AB^s +
B^sA)=\{1,e^{i2\pi/3},e^{i4\pi/3}\}$.

{\bf Subcase 6.} $c_1=0,c_2\neq0,c_3=0$.

Pick $\alpha\neq0$ such that $\alpha^2c_2+1=0$. Since the
matrix$$C=\left(\begin{array}{ccc}
\alpha&1&0\\
0&-\alpha&-1\\
0&0&0\end{array}\right)$$
has distinct eigenvalues $\alpha,-\alpha,0$, there exists an operator $B$
of rank 2 such that the operator matrix of $B^s$ equals $C\oplus0$.
Then $AB^s + B^sA$ has operator matrix $$\left(\begin{array}{cccc}
1&0&c_2&*\\
0&0&-\alpha c_2&*\\
0&-\alpha&-1&0\\
0&0&0&0\end{array}\right),$$
and $\sigma_\pi(AB^s + B^sA)=\{1,e^{i2\pi/3},e^{i4\pi/3}\}$.

{\bf Subcase 7.} $c_1=0,c_2=0,c_3=0$.

In this subcase we take  $$C=\left(\begin{array}{ccc}
-1&1&-1\\
2&1&-1\\
0&0&0\end{array}\right)$$ which has distinct eigenvalues
$\sqrt{3},-\sqrt{3},0$. Thus there exists an operator $B$ of rank 2
such that the operator matrix of $B^s$ equals $C\oplus0$. Then $AB^s
+ B^sA$ has operator matrix $$\left(\begin{array}{cccc}
1&-1&0&*\\
0&0&-1&*\\
2&1&-1&0\\
0&0&0&0\end{array}\right)$$ with $\sigma_\pi(AB^s +
B^sA)=\{3^{1/3},3^{1/3}e^{i2\pi/3},3^{1/3}e^{i4\pi/3}\}$.

{\bf Case 2.} For every $x\in X$, $\{x, Ax,A^2x\}$ is a linearly
dependent set.

In this case $A$ is a locally algebraic operator, and hence a result
due to Kaplansky  (see, for example, \cite{K2}) tells us that $A$ is
an algebraic operator of degree not greater than 2. So there exist
$\alpha,\beta,\gamma\in \mathbb{C}$ such that $\alpha A^2+\beta
A+\gamma I=0$ with $(\alpha,\beta,\gamma)\neq(0,0,0)$.

If $\alpha=0$, then $\beta\neq0$ and $\gamma\neq0$, and therefore,
$A$ is a scalar operator. Take
$B=\rm{diag}(1,e^{i2\pi/3s},e^{i4\pi/3s})\oplus0$. Then
$\sigma_\pi(AB^s + B^sA)$ contains three different points, a
contradiction.

Now assume that $\alpha\neq0$. Since $A^2\neq0$, it follows that
$(\beta,\gamma)\neq(0,0)$ and $\sigma (A)=\{a,b\}$, where
$a,b\in\mathbb{C}$.

{\bf Subcase 1.} $\sigma (A)=\{a,b\}$ with $a\neq0$ and $b\neq0$.

Then there exist linearly independent vectors $x_1,x_2\in X$ such that $Ax_1=ax_1,Ax_2=bx_2$.
Since $A$ has rank at least 3, there is $y\in X$ such that $y\notin [x_1,x_2]$ and $Ay=ay$ or $Ay=by$.
Then there is a decomposition of $X$ so that $A$ has operator matrix
$$\left(\begin{array}{cc}
A_0&*\\
0&*\end{array}\right)$$
where $A_0=\rm{diag}(a,b,a)$ or $\rm{diag}(a,b,b)$. Then there is $B$ with operator matrix $B_1\oplus0$
, where $B_1=\rm{diag}(1,b_1,b_2)$, such that $AB^s+B^sA$ has operator matrix
$$\left(\begin{array}{cc}
A_0B_1^s+B_1^sA_0&*\\
0&0\end{array}\right)$$
which has rank 3 and $\sigma_\pi(AB^s + B^sA)$ contains three different points.

{\bf Subcase 2.} $\sigma(A)=\{a,0\}$ with $a\neq0$.

Then there exist linearly independent vectors $x_1,x_2\in X$ such that $Ax_1=ax_1,Ax_2=0$.
Since $A$ has rank at least 3, there is $y\in X$ such that $y\notin [x_1,x_2]$ and $Ay=ay$.
Then there is a decomposition of $X$ so that $A$ has operator matrix
$$\left(\begin{array}{cc}
A_0&*\\
0&*\end{array}\right)$$
where $A_0=\rm{diag}(a,0,a)$.
In this case, we can use the argument in the proof when $A$ has rank 2 to choose $B$ with
operator matrix $B_1\oplus0$ so that $B_1\in M_3$ and $AB^s+B^sA$ has operator matrix
$$\left(\begin{array}{cc}
A_0B_1^s+B_1^sA_0&*\\
0&0\end{array}\right)$$
which has rank 3 and $\sigma_\pi(AB^s + B^sA)$ contains three different points.\hfill$\Box$

{\bf Corollary 2.3.} {\it  Suppose $s$ is a positive integer. Let
$X$ be  a complex Banach space with $\dim X\geq 3$ and let
$0\not=A\in\mathcal{B}(X)$. Then the following conditions are
equivalent.}

(1){ \it $A$ has rank one, or $A$ has rank two with $A^2=0$.}

(2) {\it $\sigma_\pi(AB^s+B^sA)$ has at most two elements for any $B$ in $B(X)$.}

(3){ \it There  exists no operator $B$ with rank at most 3 such that
$B^rAB^s+B^sAB^r$ has rank at most 6 and
$\sigma_\pi(B^rAB^s+B^sAB^r)$ has more than two elements.}

{\bf Proof}. (1)$\Rightarrow$(2). If $A$ has rank one, then (2) clearly holds. If $A$ has
rank two and $A^2=0$, then there is a decomposition of $X$ such that $A$ has operator matrix
$$\left(\begin{array}{ccc} 0_2&I_2&0\\0_2&0_2&0\\0&0&0\end{array}\right).$$
So, for any $B$ in $\mathcal{A}$ such that $B^s$ has operator matrix
$$\left(\begin{array}{ccc}
B_{11}&B_{12}&B_{13}\\
B_{21}&B_{22}&B_{23}\\
B_{31}&B_{32}&B_{33}\end{array}\right),$$
$AB^s+B^sA$ has operator matrix
$$\left(\begin{array}{ccc}
B_{21}&B_{11}+B_{22}&B_{23}\\
0&B_{21}&0\\
0&B_{31}&0\end{array}\right).$$ Then
$\sigma_\pi(AB^s+B^sA)=\sigma_\pi(B_{21})$ has at most two different
elements as $B_{21}\in M_2$.

The implication (2)$\Rightarrow$(3) is clear.

Finally, we verify the implication (3)$\Rightarrow$(1). If (3) holds, by Lemma 2.2,
we see that $A$ is rank 1 whenever $A^2\neq0$. If $A^2=0$, we claim that $\rm{rank}A\leq2$.
If it is not true, then we can find linearly independent vectors $x_1,x_2,x_3\in X$ such
that $\{Ax_1,Ax_2,Ax_3\}$ is a linearly independent set. It follows from $A^2=0$ that
$\{x_1,x_2,x_3,Ax_1,Ax_2,Ax_3\}$ is a linearly independent set. Let $N$ be a closed subspace
of $X$ such that $X=\rm{span}\{x_1,x_2,x_3,Ax_1,Ax_2,Ax_3\}\oplus N$.
Then $A$ has operator matrix
$$\left(\begin{array}{ccc} 0_3&0_3&*\\I_3&0_3&*\\0&0&*\end{array}\right).$$
Take $$B=\left(\begin{array}{cc} D&D\\0_3&0_3\end{array}\right)\oplus0, \mbox {with}\
D=\rm{diag}(1,e^{i2\pi/3s},e^{i4\pi/3s}).$$
Then $$\sigma_\pi(AB^s+B^sA)=\left(\begin{array}{ccc} C&0_3&*\\C&C&0\\0&0&0\end{array}\right),
\mbox {with}\ C=\rm{diag}(1,e^{i2\pi/3},e^{i4\pi/3}),$$
which has rank 6 and $\sigma_\pi(AB^s+B^sA)=\{1,e^{i2\pi/3},e^{i4\pi/3}\}$. \hfill$\Box$

The following lemma comes from \cite{CL}.

 {\bf Lemma 2.4.} {\it Let $x\in X$ and $f\in X^*$. Then, for every
$B\in \mathcal{A}$,

$$
\sigma_\pi(Bx\otimes f+x\otimes fB) = \left\{ \begin{array}{ll}
\{f(Bx)\} &{\rm if}\ \mbox\ f(x)=0 \ \mbox {\rm or}\ f(B^2x)=0,\\
\{\pm\sqrt{f(B^2x)f(x)}\} &{\rm if}\ \mbox\ f(x)\neq0,f(Bx)=0,f(B^2x)\neq0,\\
\{\alpha\} &{\rm if} \mbox\ f(x)\neq0,f(Bx)\neq0,f(B^2x)\neq0,
\end{array} \right.$$
where the scalar
$$|\alpha|=\max\{|f(Bx)+\sqrt{f(B^2x)f(x)}|,|f(Bx)-\sqrt{f(B^2x)f(x)}|\}.$$}

\section{Proof of the main result}

In this section we will complete the proof of Theorem 1.2.

It is clear that Theorem 1.2   follows from the special case below,
by considering $A_{i_p}=A$ and all other $A_{i_q}=B$.

{\bf Theorem 3.1.} {\it Let $\mathcal{A}_1$ and $\mathcal{A}_2$ be
standard operator algebras on complex Banach spaces $X_1$ and $X_2$,
respectively. Assume that
$\Phi:\mathcal{A}_1\rightarrow\mathcal{A}_2$ is a map, the range of
which contains all operators of rank at most three. Suppose also
that $\Phi$ satisfies
$$\sigma_\pi(B^rAB^s+B^sAB^r)=\sigma_\pi(\Phi(B)^r\Phi(A)\Phi(B)^s+\Phi(B)^s\Phi(A)\Phi(B)^r)
\eqno(3.1)$$   for all $ A,B\in\mathcal{A}_1$. Then one of the two
assertions in Theorem $1.2$ holds with $m=r+s+1$.}

 Thus we focus our attention to prove Theorem 3.1.

We note that the case when $s=r>0$ has been verified in \cite{ZH2}.
So, unless specified otherwise, we will assume $s>r\geq0$ in the
rest of this section. In below, we first show that $\Phi$ in Theorem
3.1 is injective.

For a Banach space $X$, denote by $\mathcal{I}_1(X)$ the set of all
rank one idempotent operators in $\mathcal{B}(X)$. In other words,
$\mathcal{I}_1(X)$ consists of all bounded operators $x\otimes f$
with $x\in X,f\in X^*$ and $\langle x,f\rangle=f(x)=1$.

The following Lemma 3.2 was proved in \cite{HLW2}.

{\bf Lemma 3.2.} {\it Let $A,A'\in \mathcal{B}(X)$ for some Banach space $X$. Suppose
$$\langle Ax,f\rangle=0\Leftrightarrow\langle A'x,f\rangle=0,\quad\mbox {\it for all}\
 x\otimes f\in\mathcal{I}_1(X).$$
Then $A'=\lambda A$ for some scalar $\lambda$.}

{\bf Lemma 3.3.} {\it Suppose $r$ and $s$ are nonnegative integers
with $(r, s)\neq(0, 0)$.
 Let $X$ be a complex Banach space. If $A,A'\in\mathcal{B}(X)$ satisfy
 $\sigma_\pi(B^rAB^s+B^sAB^r)=\sigma_\pi(B^rA'B^s+B^sA'B^r)$for all $
 B\in\mathcal{I}_1(X)$,
then $A=A'$.}

{\bf Proof.} We may suppose that $A'\neq0$ since it is obvious that
$\sigma_\pi(B^rAB^s + B^sAB^r) = \{0\}$ for all rank one idempotents $B$ implies that $A = 0$.

Assume first that $s\geq r>0$. Then the assumption implies that $\sigma_\pi(BAB)=\sigma_\pi(BA'B)$
and hence $f(Ax) = \tr(BAB) = \tr(BA'B) = f(A'x)$ for all rank one idempotents $B = x\otimes f$.
By Lemma 3.2, we see that $A' = A$.

Assume then that $s > r = 0$ and write the rank-one idempotent $B$
in the form $B = x\otimes f$ with $\langle x,f\rangle=1$. Then
$AB^s+B^sA=AB+BA$, and $\tr(AB+BA)= 0$ if and only if
$\sigma_\pi(AB+BA)=\{0\}$ or $\{\beta,-\beta\}$ for some nonzero
$\beta$. Since $\sigma_\pi(AB+BA)=\sigma_\pi(A'B+BA')$, we see that
$\tr(AB+BA)=0$ if and only if $\tr(A'B+BA')=0$. It follows from
Lemma 3.2 again that $A'=\lambda A$ for some scalar $\lambda$. But
the peripheral spectrum coincidence implies $\lambda =
1$.\hfill$\Box$

As a direct consequence of Lemma 3.3 and the condition (3.1), we have

{\bf Corollary 3.4.} {\it Let $\Phi$ satisfy the hypothesis of Theorem 3.1.
Then $\Phi$ is injective, and $\Phi(0)=0$.}

To complete the proof of Theorem 3.1, we need some more technical
lemmas.

{\bf Lemma 3.5.} {\it Let $P,Q\in\mathcal{I}_1(X)$. Then $PQ=0=QP$ if and only if
there is $B\in\mathcal{B}(X)$, which can be chosen to have rank 2, such that $\sigma_\pi
(PB+BP)=\{2\}$, $\sigma_\pi(QB+BQ)=\{-2\}$, and $\sigma_\pi(BR+RB)=\{0\}$ whenever
$R\in\mathcal{I}_1(X)$ satisfies $\sigma_\pi(PR+RP)=\sigma_\pi(QR+RQ)=\{0\}$.}

{\bf Proof.} Suppose $P,Q\in\mathcal{I}_1(X)$ satisfy $PQ=0=QP$. Then there is a space decomposition
for $X$ such that $P$ and $Q$ have operator matrices
$$\left(\begin{array}{cc} 1&0\\0&0\end{array}\right)\oplus0 \quad\mbox {\rm and}\quad\mbox\
\left(\begin{array}{cc} 0&0\\0&1\end{array}\right)\oplus0.$$
Using the same space decomposition, let $B$ have operator matrix
$\left(\begin{array}{cc} 1&0\\0&-1\end{array}\right)\oplus0$. Then $B$ has rank 2 such that
$\sigma_\pi(PB+BP)=\{2\}$ and $\sigma_\pi(QB+BQ)=\{-2\}$. Consider any $R\in\mathcal{I}_1(X)$
such that $\sigma_\pi(PR+RP)=\sigma_\pi(QR+RQ)=\{0\}.$ Using the same space decomposition
as $P$ and $Q$, we assume that $R$ has operator matrix
$$\left(\begin{array}{cc} R_{11}&R_{12}\\R_{21}&R_{22}\end{array}\right)$$
where $R_{11}\in M_2$. Since
$\sigma_\pi(PR+RP)=\sigma_\pi(QR+RQ)=\{0\}$, the $(1, 1)$ and $(2,
2)$ entry of $R_{11}$ are both zero. Thus, $R_{22}$ has trace one
and rank one. We may then assume that $R_{22}$ has operator matrix
$(1)\oplus0$. As a result, we may assume that the operator matrix of
$R$ has the form $\hat{R}\oplus0$, where $R$ or $R^t$ has one of the
following forms:
$$\left(\begin{array}{ccc}0&0&*\\0&0&*\\0&0&1\end{array}\right)\quad\mbox {\rm or}\quad\mbox\
\left(\begin{array}{ccc}0&a&b\\0&0&0\\0&c&1\end{array}\right)\quad\mbox {\rm with}\quad\mbox\
a=bc.$$
Consequently, $\sigma_\pi(BR+RB)=\{0\}$.

Conversely, suppose $P,Q\in\mathcal{I}_1(X)$ such that $PQ\neq0$ or
$QP\neq0$. Then there is a space decomposition for $X=X_1\oplus X_2$
with $\dim X_1=2$ such that $P$ has operator matrix
$\left(\begin{array}{cc} 1&0\\0&0\end{array}\right)\oplus0$ and $Q$
has operator matrix
$$\left(\begin{array}{cc} 0&0\\1&1\end{array}\right)\oplus0 \quad\mbox {\rm or}\quad\mbox\
\left(\begin{array}{cc} 0&1\\0&1\end{array}\right)\oplus0.$$ We
assume that the former case holds. The proof for the other case is
similar. Suppose there is a $B$ in $\mathcal{B}(X)$ such that
$\sigma_\pi (PB+BP)=\{2\}$, $\sigma_\pi(QB+BQ)=\{-2\}$ and
$\sigma_\pi(BR+RB)=\{0\}$ whenever $R$ in $\mathcal{I}_1(X)$
satisfies $\sigma_\pi(PR+RP)=\sigma_\pi(QR+RQ)=\{0\}$. Using the
same space decomposition as $P$ and $Q$, we assume that $B$ has
operator matrix
$$\left(\begin{array}{cc} B_{11}&B_{12}\\B_{21}&B_{22}\end{array}\right)$$
where $B_{11}\in M_2$.

First, we claim that $B_{22}=0$. If not, we may assume that the $(1,1)$ entry of $B_{22}$
is nonzero. If $R$ has operator matrix
$$\left(\begin{array}{ccc} 0&0&0\\0&0&0\\0&0&1\end{array}\right)\oplus0,$$
we see that $\sigma_\pi(PR+RP)=\sigma_\pi(QR+RQ)=\{0\}\neq\sigma_\pi(BR+RB).$

 Next, we claim that $B_{12}=0$. If this is not true, we can find a
suitable space decomposition for $X_2$ such that $B_{12}$ has the
form $\left(\begin{array}{cc} 1&0\\0&T\end{array}\right),$ where the
last column is vacuous if dim$X=3$, and $T$ has rank zero or one.
But then if $R\in\mathcal{I}_1(X)$ has operator matrix
$$\left(\begin{array}{ccc} 0&0&0\\0&0&0\\1&0&1\end{array}\right)\oplus0,$$
we have $\sigma_\pi(BR+RB)\neq\{0\}.$ Similarly, we can show that $B_{21} = 0$.

Now, we consider $B_{11}=\left(\begin{array}{cc}
b_{11}&b_{12}\\b_{21}&b_{22}\end{array}\right).$ Let $R$ has
operator matrix $\left(\begin{array}{ccc}
0&0&0\\1&0&1\\1&0&1\end{array}\right)\oplus0,$ We see that
$\sigma_\pi(PR+RP)=\sigma_\pi(QR+RQ)=\{0\}$. Because $BR+RB$ has
operator matrix
$$\left(\begin{array}{ccc}
b_{12}&0&b_{12}\\
b_{11}+b_{22}&b_{12}&b_{22}\\
b_{11}&b_{12}&0\end{array}\right)\oplus0,$$
so $b_{12}=0$. Since $\sigma_\pi(PB+BP)=\{2\}$ and $\sigma_\pi(QB+QB)=\{-2\}$,
it follows that $b_{11}=1,b_{22}=-1$. Finally, for $R$ with operator
matrix $\left(\begin{array}{cc} 0&-1\\0&1\end{array}\right)\oplus0,$ we have
$\sigma_\pi(PR+RP)=\sigma_\pi(QR+RQ)=\{0\}$. But $BR+RB$ has operator matrix
$$\left(\begin{array}{cc} -b_{21}&0\\b_{21}&-2-b_{21}\end{array}\right)\oplus0,$$
which cannot be a nilpotent.\hfill$\Box$

For a Banach space $X$ and a ring automorphism $\tau$ of $\mathbb{C}$, if an additive map
$T : X\rightarrow X$ satisfies $T(\lambda x) = \tau(\lambda)Tx$ for all complex $\lambda$
and all vectors $x$, we say that $T$ is $\tau$-linear.

The following result can be proved by a similar argument to the
proof of the main result in \cite{M2}, see also \cite{BH} and
\cite{S2}.

{\bf Lemma 3.6.} {\it Let $X_1$ and $X_2$ be complex Banach spaces
with dimension at least 3. Let
$\Phi:\mathcal{I}_1(X_1)\rightarrow\mathcal{I}_1(X_2)$ be a
bijective map with the property that
$PQ=QP=0\Leftrightarrow\Phi(P)\Phi(Q)=\Phi(Q)\Phi(P)=0$ for all
$P,Q\in\mathcal{I}_1(X_1)$. Then there exists a ring automorphism
$\tau$ of $\mathbb{C}$ such that one of the following cases holds.}

(i) {\it There exists a $\tau$-linear transformation $T :
X_1\rightarrow X_2$ satisfying $\Phi(P)= TPT^{-1}$ for all
$P\in\mathcal{I}_1(X_1)$.}

(ii) {\it There exists a $\tau$-linear transformation $T :
X^*_1\rightarrow X_2$ satisfying $\Phi(P)=TP^*T^{-1}$ for  all
$P\in\mathcal{I}_1(X_1)$.

If $X$ is infinite dimensional, the transformation $T$ is an
invertible bounded linear or conjugate linear operator.}\vskip 4mm

In the following, we present the proof of Theorem 3.1.

{\bf Proof of Theorem 3.1.} Recall that $\Phi$ satisfies condition
(3.1).

{\bf Case 1.} $s>r>0$.

{\bf Claim 1.1.} $\Phi$ is injective, and $\Phi(0)=0$.

It is just Corollary 3.4.

{\bf Claim 1.2.} $\Phi$ preserves rank one operators in both
directions.

Assume that rank$A=1$; then Claim 1.1 implies that $\Phi(A)\neq0$.
For any $B\in {\mathcal A}_1$, by Lemma 2.1,
$\sigma_\pi(\Phi(B)^r\Phi(A)\Phi(B)^s+\Phi(B)^s\Phi(A)\Phi(B)^r)=\sigma_\pi(B^rAB^s+B^sAB^r)$
has at most two different elements. Since the range of $\Phi$
contains all operators of rank at most three, for any $C\in{\mathcal
A}_2$ with rank$(C)\leq 3$,
$\sigma_\pi(C^r\Phi(A)C^s)+C^s\Phi(A)C^r)$ has at most two different
elements. Applying Lemma 2.1 again one sees that $\Phi(A)$ is of
rank one. Conversely, assume that  $T\in{\mathcal A}_2$ is of rank
one. Then, there is $A\in{\mathcal A}_1$ such that $\Phi(A)=T$. For
any $B\in{\mathcal A}_1$, Lemma 2.1 implies that
$\sigma_\pi(B^rAB^s+B^sAB^r)=\sigma_\pi(\Phi(B)^r\Phi(A)\Phi(B)^s)+\Phi(B)^s\Phi(A)\Phi(B)^r)$
has at most two different elements. Applying Lemma 2.1 again one
gets $A$ is of rank one.

{\bf Claim 1.3.} $\Phi$ is linear.

We show first that $\Phi$ is additive.

Observe that, for any operator $A$ and rank one operator $x\otimes
f$, we have
$$(x\otimes f)^rA(x\otimes f)^s+(x\otimes f)^sA(x\otimes f)^r=2\langle
x,f\rangle^{s+r-2}(x\otimes f) A(x\otimes f),
$$
and hence
$$\begin{array}{rl} &\sigma_\pi((x\otimes f)^rA(x\otimes f)^s+(x\otimes f)^sA(x\otimes
f)^r) \\= &{\rm Tr}(2\langle x,f\rangle^{s+r-2}(x\otimes f)
A(x\otimes f))=\{2\langle x,f\rangle^{s+r-1}\langle
Ax,f\rangle\}.\end{array} \eqno(3.3)$$

Let $A,B\in \mathcal{A}_1$ be arbitrary. For any $y\in X_2, g\in
X_2^*$ with $\langle y,g\rangle=1$, Claim 2 implies that there exist
$x\in X_1, f\in X_1^*$ such that $\Phi(x\otimes f)=y\otimes g$.

Then, by Eqs.(3.1) and (3.3), we have
$$\begin{array}{rl}
\{2\langle\Phi(A+B)y, g\rangle\langle y, g\rangle\}
=&\sigma_\pi((y\otimes g)^r\Phi(A+B)(y\otimes g)^s+(y\otimes g)^s\Phi(A+B)(y\otimes g)^r)\\
=&\sigma_\pi((x\otimes f)^r(A+B)(x\otimes f)^s+(x\otimes f)^s(A+B)(x\otimes f)^r)\\
=&\{{\rm Tr}(2\langle x,f\rangle^{r+s-2}(x\otimes f)(A+B)(x\otimes f))\}\\
=&\{{\rm Tr}(2\langle x,f\rangle^{r+s-2}(x\otimes f)A(x\otimes f))\} \\
&+\{{\rm Tr}(2\langle x,f\rangle^{r+s-2}(x\otimes f)B(x\otimes f))\}\\
=&\{{\rm Tr}(2\langle y,g\rangle^{s+r-2}(y\otimes g)  \Phi(A) (y\otimes g)) \\
&+\{{\rm Tr}(2\langle y,g\rangle^{s+r-2}(y\otimes g)  \Phi(B) (y\otimes g))\\
=& \{{\rm Tr}(2\langle y,g\rangle^{s+r-2}(y\otimes g) (\Phi(A)+\Phi(B))(y\otimes g))\\
=&\{\langle (\Phi(A)+\Phi(B))y, g\rangle\langle y,
g\rangle\}.\end{array}$$  It follows that
$$\langle\Phi(A+B)y,g\rangle=\langle(\Phi(A)+\Phi(B))y, g\rangle$$
holds for any $y\in X_2, g\in X_2^*$ with $\langle y,g\rangle=1$.
This entails $\Phi(A+B)=\Phi(A)+\Phi(B)$ and hence $\Phi$ is
additive. Similarly one can check that $\Phi$ is homogeneous, I.e.,
$\Phi(\lambda A)=\lambda\Phi(A)$. So $\Phi$ is linear.

The claims 1.1-1.3 imply that  $\Phi$ is an injective linear map
preserving rank one operators in both directions. By \cite{H1} the
following claim is true.

{\bf Claim 1.4.} One of the following statements holds:

(i) There exist two linear bijections $T:X_1\rightarrow X_2$ and
$S:X^*_1\rightarrow X^*_2$ such that $\Phi(x\otimes f)=Tx\otimes Sf$ for
all rank one operators $x\otimes f\in \mathcal{A}_1$.

(ii) There exist two linear bijections $T:X^*_1\rightarrow X_2$ and
$S:X_1\rightarrow X^*_2$ such that $\Phi(x\otimes f)=Tf\otimes Sx$ for
all rank one operators $x\otimes f\in \mathcal{A}_1$.

{\bf Claim 1.5.} There exists a scalar $\lambda\in \mathbb{C}$ with
$\lambda^m=1$ and $m=r+s+1$ such that, if (i) occurs in Claim 1.4,
then $\langle Tx, Sf\rangle =\lambda \langle x,f\rangle$ holds for
all $x\in X_1$ and $f\in X^*_1$; if (ii) occurs in Claim $4$, then
$\langle Tf,Sx\rangle =\lambda \langle x,f\rangle$ holds for all
$x\in X_1$ and $f\in X^*_1$.

To check Claim 1.5, we first assume that the case (i) in Claim 1.4
occurs. Then, for any $x\in X_1$, $ f\in X^*_1$, we have
$\sigma_\pi(2(x\otimes f)^m)=\{2{\langle x, f\rangle}^m\}
=\sigma_\pi(2(Tx\otimes Sf)^m)=\{2{\langle Tx, Sf\rangle}^m\}$. So
$\langle Tx, Sf\rangle=\lambda_{x,f}\langle x, f\rangle$ with
${\lambda_{x,f}}^m=1$. Especially, $\langle x,
f\rangle=0\Leftrightarrow\langle Tx, Sf\rangle=0$.

Let $V_0=\{(x,f)\mid \langle x,f\rangle=0\}$, $V_t=\{(x,f)\mid\lambda_{x,f}
=e^{i\frac{2(t-1)\pi}{m}}\}$, $t=1,\ldots,m$.
 Then $\bigcup_{t=1}^m V_t=X\times X^*$ and $V_k\cap V_j=V_0$ if $k\not=j$,
$k,j=1,2,\ldots,m$. For $x_1,x_2\in X_1$ we denote by $[x_1,x_2]$ the linear
subspace spanned by $x_1$ and $x_2$.

 {\bf Assertion 1.} For any nonzero $x_1,x_2\in X_1, f\in X^*_1$,
 there exists some $k\in\{1,2,\ldots,m\}$ such that
$[x_1,x_2]\times[f]\subseteq V_k$.

We need only to show that we may take $\lambda_{x_1,f}$ and
$\lambda_{x_2,f}$ such that $\lambda_{x_1,f}=\lambda_{x_2,f}$. Consider the following
three cases.

{\bf Case 1$^\circ$.} $x_1,x_2$ are linearly dependent.

Assume that $x_2=\alpha x_1$; then $\alpha\not=0$ and
$\alpha\lambda_{x_1,f}\langle x_1, f\rangle=\alpha\langle Tx_1,
Sf\rangle=\langle Tx_2, Sf\rangle=\alpha\lambda_{x_2,f}\langle x_1,
f\rangle$. So we may take $\lambda_{x_1,f}$ and $\lambda _{x_2,f}$ such that
$\lambda_{x_1,f}=\lambda _{x_2,f}$.

 {\bf Case 2$^\circ$.} $x_1,x_2$ are
linearly independent and at least one of $\langle x_i, f\rangle$,
$i=1,2$ is not   zero.

In this case, for any $\alpha,\beta\in \mathbb{C}$ we have
$$\alpha\lambda_{\alpha,\beta}\langle x_1,
f\rangle+\beta\lambda_{\alpha,\beta}\langle x_2,f\rangle=\langle
T(\alpha x_1+\beta x_2), Sf\rangle=\alpha\lambda_{x_1,f}\langle
x_1,f\rangle+\beta\lambda_{x_2,f}\langle x_2,f\rangle,\eqno(3.2)$$
where $\lambda_{\alpha,\beta}=\lambda_{\alpha x_1+\beta x_2,f}.$ Let
$$\eta=\left(\begin{array}{c}
\lambda_{x_1,f}\langle x_1,f\rangle\\
\lambda_{x_2,f}\langle x_2,f\rangle
\end{array}\right),\
\eta_0=\left(\begin{array}{c}
\langle x_1,f\rangle\\
\langle x_2,f\rangle
\end{array}\right),\
\xi=\left(\begin{array}{c}
\alpha\\
\beta
\end{array}\right)\in \mathbb{C}^2.$$
Then Eq.(3.2) implies that
$$\langle\eta,\xi\rangle=\lambda_{\alpha,\beta}\langle\eta_0,\xi\rangle$$
holds for any $\xi\in\mathbb{C}^2$. It follows that
$\langle\eta,\xi\rangle=0\Leftrightarrow\langle\eta_0,\xi\rangle=0$.
So, as the vectors in ${\mathbb C}^2$, we must have $\eta=\gamma\eta_0$
for some scalar $\gamma$. Now it is clear that
$\lambda_{x_1,f}=\lambda_{x_2,f}$.

{\bf Case 3$^\circ$.}  $x_1,x_2$ are linearly independent and
$\langle x_1, f\rangle=\langle x_2,f\rangle=0$.

Then $\langle Tx_1,Sf\rangle=\lambda_{x_1,f}\langle x_1,f\rangle=0=\langle
Tx_2,Sf\rangle=\lambda_{x_2,f}\langle x_2,f\rangle$. In this case it
is clear that we can take $\lambda_{x_1,f}$ and $\lambda_{x_2,f}$
such that $\lambda_{x_1,f}=\lambda_{x_2,f}$.

Similar to the previous discussion, we have

{\bf Assertion 2.} For any nonzero $x\in X_1$, $f_1,f_2\in X^*_1$, there exists some
$k\in\{1,2,\ldots,m\}$ such that $[x]\times[f_1,f_2]\subseteq V_k$.

 {\bf Assertion 3.} There exists a scalar $\lambda\in\mathbb{C}$ with $\lambda^m=1$
such that $\lambda_{x,f}=\lambda$ for all $x\in X_1$ and $f\in X^*_1$.

For any $f_0\neq0$, there exist $x_0$ such that $\langle x_0,
f_0\rangle=1$. Then $\langle Tx_0, Sf_0\rangle=\lambda_{x_0,f_0}$
and $(x_0,f_0)\in V_{k_0}$ for some $k_0\in\{1,2,\ldots,m\}$. \if For any
$x\in\ker f_0$, that is $\langle x,f_0\rangle=0$, we may take
$\lambda_{x,f_0}=\lambda_{x_0,f_0}$. Then $(x,f_0)\in V_{k_0}$.\fi
So, by Assertion 1, for any $x\in X_1$, we have $[x,x_0]\times
[f_0]\subseteq V_{k_0}$, which implies that $X_1\times[f_0]\subseteq
V_{k_0}$. Similarly, by Assertion 2 one gets,   for any $x_0\neq0$,
$[x_0]\times X^*_1\subseteq V_{k_0}$. Thus we obtain that $X_1\times
X^*_1=V_{k_0}$.

Hence, there exists a scalar $\lambda\in \mathbb{C}$ with $\lambda^m=1$ such
that $\lambda_{x,f}=\lambda$ for all $x$ and $f$, that is, $\langle
Tx, Sf\rangle =\lambda \langle x,f\rangle$ holds for all $x\in X_1$
and $f\in X^*_1$. So Assertion 3 is true.

This completes the proof  of Claim 1.5 for the case (i) of Claim
1.4.

If the case (ii) in Claim 1.4 occurs, by a similar argument one can
show that there exists a scalar $\lambda$ with $\lambda^m=1$ such
that $\langle Tf,Sx\rangle =\lambda \langle x,f\rangle$ holds for
all $x\in X_1$ and $f\in X^*_1$. Hence the last conclusion of Claim
1.5 is also true.

{\bf Claim 1.6.} There exists a scalar $\lambda $ with $\lambda^m=1$
such that one of the followings holds:

(1) There exists an invertible operator $T\in {\mathcal B}(X_1, X_2)$
such that $\Phi(x\otimes f)=\lambda T(x\otimes f)T^{-1}$ for all
$x\otimes f\in \mathcal{A}_1$.

(2) $X_1$ and $X_2$ are reflexive, and there exists an invertible
operator $T\in {\mathcal B}(X^*_1, X_2)$ such that $\Phi(x\otimes
f)=\lambda T(x\otimes f)^* T^{-1}$ for all $x\otimes f\in
\mathcal{A}_1$.

Suppose that the case (i) of Claim 1.4 occurs. Then by Claim 1.5,
there exists a scalar $\lambda\in \mathbb{C}$ with $\lambda^m=1$
such that $\langle Tx, Sf\rangle =\lambda \langle x,f\rangle$ holds
for all $x\in X_1$ and $f\in X^*_1$. If $\{x_n\}\subset X_1$ is a
sequence such that $x_n\rightarrow x$ and $Tx_n\rightarrow y$ for
some $x\in X_1$ and $y\in X_2$ as $n\rightarrow\infty$, then, for
any $f\in X^*_1$, we have $\langle y,
Sf\rangle=\lim_{n\rightarrow\infty}\langle Tx_n, Sf\rangle=
\lim_{n\rightarrow\infty}\lambda\langle x_n, f\rangle=\lambda\langle
x,f\rangle=\langle Tx, Sf\rangle.$ As $S$ is surjective we must have
$y=Tx$. So the bijection $T$ is a closed operator and thus  a
bounded invertible operator. Since $\langle Tx, Sf\rangle=\langle x,
T^*Sf\rangle=\lambda\langle x, f\rangle$ holds for all $x\in X_1$
and $f\in X^*_1$, we see that $T^*S=\lambda I$, that is
$S=\lambda(T^*)^{-1}$. It follows from the case (i) of Claim 1.4
that $\Phi(x\otimes f)=Tx\otimes Sf=\lambda Tx\otimes
(T^*)^{-1}f=\lambda T(x\otimes f)T^{-1}$ for any rank one operator
$x\otimes f$, i.e., the case (1) of Claim 1.6 holds.

Suppose that the case (ii) of Claim 1.4 occurs. Then by Claim 1.5,
there exists a scalar $\lambda\in \mathbb{C}$ with $\lambda^m=1$
such that $\langle Tf, Sx\rangle =\lambda \langle x,f\rangle$ holds
for all $x\in X_1$ and $f\in X^*_1$. Similar to the above argument
we can check that both $T$ and $S$ are bounded invertible operators
with $S=\lambda (T^*)^{-1}$. It follows that $\Phi(x\otimes
f)=\lambda T(x\otimes f)^*T^{-1}$ for any $x\otimes f$, obtaining
that the case (2) of Claim 1.6 holds. Moreover, by \cite{H1}, in
this case both $X_1$ and $X_2$ are reflexive.

{\bf Claim 1.7.} The theorem is true.

Assume that we have the case (1) of Claim 1.6.  Let $A\in
\mathcal{A}_1$ be arbitrary. For any $x\in X_1$ and $f\in X^*_1$
with $\langle x, f\rangle=1$, we have
$$\begin{array}{rl}
\{2\langle Ax,f\rangle\}
=&\sigma_\pi((x\otimes f)^rA(x\otimes f)^s+(x\otimes f)^sA(x\otimes f)^r)\\
=&\sigma_\pi((\lambda T(x\otimes f)T^{-1})^r\Phi(A)(\lambda T(x\otimes f)T^{-1})^s\\
&+(\lambda T(x\otimes f)T^{-1})^s\Phi(A)(\lambda T(x\otimes f)T^{-1})^r)\\
=&\sigma_\pi(\frac{2}{\lambda}(x\otimes f)T^{-1}\Phi(A)T(x\otimes f))\\
=&\{\langle\frac{2}{\lambda}T^{-1}\Phi(A)Tx,f\rangle\}.\end{array}$$
This implies that $\Phi(A)=\lambda TAT^{-1}$ for any $A\in \mathcal{A}_1$.

A similar argument shows that if the case (2) of Claim 1.6 occurs
then $\Phi$ has the form given in (2) of Theorem 1.2.

{\bf Case 2.} $s>r=0$.

{\bf Claim 2.1.} $\Phi$ is injective, and $\Phi(0)=0$.

It is just Corollary 3.4.

{\bf Claim 2.2.} If $A\in\mathcal{A}_1$ is a nonzero multiple of a
rank one idempotent, then so is $\Phi(A)$. In particular, if
$P\in\mathcal{I}_1(X_1)$, then $\Phi(P)=\mu R$ such that
$R\in\mathcal{I}_1(X_2)$ and $\mu^{s+1}=1$.

Let $A\neq0$ be a nonzero multiple of an idempotent, say $A=\alpha
P$, where $0\neq\alpha\in\mathbb{C}$ and $P\in\mathcal{I}_1(X_1)$.
For any $D$ in $\mathcal{A}_2$ of rank at most 3, there is $C$ in
$\mathcal{A}_1$ such that $\Phi(C)=D$. By equation (3.1) we have
$\sigma_\pi(\Phi(A)D^s+D^s\Phi(A))=\sigma_\pi(AC^s+C^sA)$, which
contains at most two different elements. Putting $B=A$ in equation
(3.1), we have
$\sigma_\pi(2\Phi(A)^{s+1})=\sigma_\pi(2A^{s+1})\neq\{0\}$. Applying
Corollary 2.3, depending on $s>r=0$, we see that $\Phi(A)$ is a
nonzero multiple of rank one idempotent. Thus $\Phi$ preserves
nonzero multiples of rank one idempotents. If $P$ in $\mathcal{A}_1$
is a rank one idempotent, then $\Phi(P)=\mu R$, where
$R\in\mathcal{I}_1(X_2)$ and $\mu\in\mathbb{C}$. Since
$\{2\}=\sigma_\pi(2P^{s+1})=\sigma_\pi(2\Phi(P)^{s+1})=
\{2\mu^{s+1}\}$, we see that $\mu^{s+1}=1$.

{\bf Claim 2.3.} There exists a scalar $\lambda$ with
$\lambda^{s+1}=1$ such that $\lambda^{-1}\Phi$ sends rank one
idempotents to rank one idempotents.

Let $0\neq f\in X^*_1$. Assume that $\langle x_1,f\rangle=\langle
x_2,f\rangle=1$. By Claim 2, $\Phi(x_1\otimes f)=\lambda_1y_1\otimes
g_1$ and $\Phi(x_2\otimes f)=\lambda_2y_2\otimes g_2$, where
$g_1(y_1)=g_2(y_2)=1$ and $\lambda_1^{s+1}=\lambda_2^{s+1}=1$. Using
the peripheral spectrum equation (3.1) we have
$$\begin{array}{rl}
&\sigma_\pi(\lambda_1^s\lambda_2((y_1\otimes g_1)(y_2\otimes g_2)+(y_2\otimes g_2)(y_1\otimes g_1))\\
=&\sigma_\pi(\lambda_1\lambda_2^s(y_1\otimes g_1)(y_2\otimes g_2)+(y_2\otimes g_2)(y_1\otimes g_1))\\
=&\sigma_\pi((x_1\otimes f)(x_2\otimes f)+(x_2\otimes f)(x_1\otimes f))\\
=&\{2\}.\end{array}$$ Then
$\lambda_1^s\lambda_2=\lambda_2^s\lambda_1$. In particular,
$\lambda_1^2=\lambda_2^2$ as $\lambda_1^{s+1}=1=\lambda_2^{s+1}$.
Suppose $\lambda_1=-\lambda_2$, then we have $\sigma_\pi((y_1\otimes
g_1)(y_2\otimes g_2)+(y_2\otimes g_2)(y_1\otimes g_1))=\{-2\}$, but
by Lemma 2.4, this is impossible. So, $\lambda_1=\lambda_2$. Denote
this common value by $\lambda_f$. Similarly, for any nonzero $x$ in
$X_1$ we will have an $m$th root ($m=s+1$) $\lambda_x$ of unity
depending only on $x$ such that $\Phi(x\otimes f)=\lambda_x
Q_{x\otimes f}$, for some rank one idempotent $Q_{x\otimes f}$
whenever $f(x)=1$.

Now consider any two rank one idempotents $x_1\otimes f_1$ and $x_2\otimes f_2$ in $\mathcal{A}_1$.
We write $x_1\otimes f_1\sim x_2\otimes f_2$ if there is a scalar $\lambda$ with $\lambda^{s+1}=1$
such that $\lambda\Phi(x_i\otimes f_i)$ is a rank one idempotent for $i = 1, 2$. In case $\alpha
=\langle x_1,f_2\rangle\neq0$, we see that
$$x_1\otimes f_1\sim x_1\otimes \frac{f_2}{\alpha}=\frac{x_1}{\alpha}\otimes f_2\sim x_2\otimes f_2.$$
In case $\langle x_1,f_2\rangle=\langle x_2,f_1\rangle=0$, we also have
$$x_1\otimes f_1\sim(x_1+x_2)\otimes f_1\sim (x_1+x_2)\otimes f_2\sim x_2\otimes f_2.$$
So, Claim 2.3 is true.

By Claim 2.3, without loss of generality, we assume that $\Phi$
preserves rank one idempotents.

{\bf Claim 2.4.} If $\Phi(A)\in\mathcal{A}_2$ is a rank one
idempotent, then $A\in\mathcal{A}_1$
 is a rank one idempotent.

Assume that $\Phi(A)$ is a rank one idempotent. Suppose $A$ is not a
rank one idempotent, i.e., $A$ is rank one nilpotent or $A$ has rank
at least two. Putting $B=A$ in equation (3.1), we have
$\sigma_\pi(2A^{s+1})=\sigma_\pi(2\Phi(A)^{s+1})=\{2\}$, so
$A^2\neq0$. Then $A$ has rank at least 2. In this case, by the
arguments in the proof of Lemma 2.2, for such $A$ we can find a
operator $B$ with rank$B\leq3$ such that $\sigma_\pi(AB^s+B^sA)$ has
three different points. But
$\sigma_\pi(\Phi(A)\Phi(B)^s+\Phi(B)^s\Phi(A))$ has at most two
different points, a contradiction.

{\bf Claim 2.5.} One of the following statements is true.

(i) There exists a bounded invertible linear operator $T : X_1\rightarrow X_2$ such that
$$\Phi(x\otimes f)=T(x\otimes f)T^{-1}\quad\mbox {\rm for all } \mbox\ x\in X_1,
f\in X^*_1 \ \rm{ with} \ \langle x,f\rangle=1.$$

(ii) There exists a bounded invertible linear operator $T : X^*_1\rightarrow X_2$ such that
$$\Phi(x\otimes f)=T(x\otimes f)^*T^{-1}\quad\mbox {\rm for all }  \mbox\ x\in X_1,
f\in X^*_1 \ \rm{with}\ \langle x,f\rangle=1.$$

Since $\Phi$ preserves rank one idempotents in both directions, by use of Lemma 3.5, it is easily
checked that $P,Q\in\mathcal{I}_1(X_1)$ satisfy $PQ=0=QP$ if and only if
$\Phi(P)\Phi(Q)=0= \Phi(Q)\Phi(P)$.
Thus we can apply Lemma 3.6 to conclude that (i) or (ii) holds, but with $T$ a $\tau$-linear for
some ring automorphism $\tau$ of $\mathbb{C}$.

Next we prove that $\tau$ is the identity and hence $T$ is linear.
For any $\alpha\in \mathbb{C}\backslash\{0,1\}$, let $A$ and $B$
have respectively operator matrices
$$\left(\begin{array}{cc}
1 & \alpha-1\\
0 &  0\end{array}\right)\oplus0\quad\mbox {\rm and }\quad\mbox\
\left(\begin{array}{cc}
1 & 0\\
1 & 0\end{array}\right)\oplus0.$$ Then
$$AB^s+B^sA=\left(\begin{array}{cc}
\alpha+1 & \alpha-1\\
1 &  \alpha-1\end{array}\right)\oplus0.$$
Since $$\begin{array}{rl}
&\sigma_\pi(AB^s+B^sA)=\sigma_\pi(\Phi(A)\Phi(B)^s+\Phi(B)^s\Phi(A))\\
=&\sigma_\pi(T(AB^s+B^sA)T^{-1})=\{\tau(\xi):\xi\in\sigma_\pi(AB+BA)\},
\end{array}$$

Note that $\sigma_\pi(AB^s+B^sA)=\{\alpha\pm\sqrt{\alpha}\}
=\sigma_\pi(\tau(\alpha\pm\sqrt{\alpha}))$. It follows that
$\tau(2\alpha)=\tau(\alpha+\sqrt{\alpha}+\alpha-\sqrt{\alpha})=\alpha+\sqrt{\alpha}+\alpha-\sqrt{\alpha}=
2\alpha$. Hence $\tau(\alpha)=\alpha$ for any $\alpha\in\mathbb{C}$.
It follows that $T$ is an invertible bounded linear operator.

\if If $\alpha<0$, then
$\sigma_\pi(AB^s+B^sA)=\{\alpha\pm\sqrt{-\alpha}i\}
=\sigma_\pi(\tau(\alpha+\sqrt{-\alpha}i),\tau(\alpha-\sqrt{-\alpha}i))$.
It follows that
$\tau(2\alpha)=\tau(\alpha+\sqrt{-\alpha}i)+\tau(\alpha-\sqrt{-\alpha}i)=2\alpha$.
Thus $\tau(\alpha)=\alpha$.

If $\alpha>0$, $\tau(\alpha)=\tau(-1)\tau(-\alpha)=\alpha$.

If $\alpha=i$,
$\sigma_\pi(AB+BA)=\{\frac{1}{2}+\frac{3i}{2}\}=\{\tau(\frac{1}{2}+\frac{3i}{2})\}$,
it follows from
$\tau(\frac{1}{2}+\frac{3i}{2})=\frac{1}{2}+\frac{3i}{2}=\frac{1}{2}+\frac{3}{2}\tau(i)$
that $\tau(i)=i$.\fi

{\bf Claim 2.6.} $\Phi$ has the form in Theorem 3.1.

Suppose (i) in Claim 2.5 holds. Let $A\in\mathcal{A}_1$ be
arbitrary. For any $x\in X_1$ and $f\in X^*_1$ with $\langle
x,f\rangle=1$, the condition (3.1) ensures that
$$\begin{array}{rl}
&\sigma_\pi(T^{-1}\Phi(A)T(x\otimes f)^s+(x\otimes f)^sT^{-1}\Phi(A)T)\\
=&\sigma_\pi(T[T^{-1}\Phi(A)T(x\otimes f)^s+(x\otimes f)^sT^{-1}\Phi(A)T]T^{-1})\\
=&\sigma_\pi(\Phi(A)T(x\otimes f)^sT^{-1}+T(x\otimes f)^sT^{-1}\Phi(A))\\
=&\sigma_\pi(A(x\otimes f)^s+(x\otimes f)^sA)
\end{array},$$
Hence, by Lemma 3.3, we have $T^{-1}\Phi(A)T=A$ for all $A$ in
$\mathcal{A}_1$, that is, $\Phi$ has the form (1) in the Theorem
1.2.

Similarly, one can show that $\Phi$ has the form (2) if (ii) of
Claim 2.5 holds, completing the proof of Theorem 3.1. \hfill$\Box$

\end{document}